\documentclass
[12pt]{article}
\usepackage{epsfig}

\usepackage{amssymb,amsmath,float}
\begin{document}
\title{ \bf On Fuzzy $\Gamma$-Hypersemigroups }
\author {\bf R. Ameri$^\dag$, R. Sadeghi$^\ddag$\\
         \normalsize\sl
\small  $^\dag$ School of Mathematics, Statistic and Computer Sciences, College of Science,\\
\normalsize\sl\small University of Tehran, P.O. Box 14155-6455,
Tehran, Iran, E-mail: rameri@ut.ac.ir, \normalsize\sl\\
\normalsize\sl \small  $^\ddag$ Department of Mathematics, Aliabad
Katoul Branch, Islamic Azad University,\\ \sl\small Aliabad Katoul,
Iran, \small e-mail: $razi_{-}sadeghi@yahoo.com$}
%\singlespace

\date{}
\maketitle

\begin{abstract}

\noindent We introduced and  study fuzzy $\Gamma$-hypersemigroups,
according to fuzzy semihypergroups as previously defined \cite{S1}
and prove that results in this respect.  In this regard first we
introduce  fuzzy hyperoperation and then study fuzzy
$\Gamma$-hypersemigroup.  We will proceed by study fuzzy
$\Gamma$-hyperideals and fuzzy $\Gamma$-bihyperideals. Also we study
the relation between the classes of fuzzy $\Gamma$-hypersemigroups
and $\Gamma$-semigroups. Precisely, we  associate a
$\Gamma$-hypersemigroup to every fuzzy $\Gamma$-hypersemigroup and
vice versa. Finally, we introduce and study fuzzy
$\Gamma$-hypersemigroups regular and fuzzy strongly regular
relations of fuzzy $\Gamma$-hypersemigroups.
\end{abstract}

\noindent Keywords: fuzzy $\Gamma$-hyperoperation,
$\Gamma$-hypersemigroup, fuzzy $\Gamma$-hypersemigroup, fuzzy
$\Gamma$-hyperideals, fuzzy regular relation

\section{Introduction}
\noindent Hyperstructure theory was born in 1934 when Marty
\cite{M1} defined hypergroups, began to analysis their properties
and applied them to groups. Algebraic hyperstructures are a suitable
generalization of classical algebraic structures. In 1986, M.K. Sen
and Saha \cite{S2} defined the notion of a $\Gamma$-semigroup as a
generalization of a semigroup. Ameri \cite{A12} introduced and study
fuzzy ideals of gamma-hyperrings, after that Davvaz et. al.
\cite{D4} studied $\Gamma$-semihypergroup as a generalization of a
semihypergroup and then many classical notions of semigroups and
semihypergroups have been extended to $\Gamma$-semihypergroups.
Zadeh \cite{Z1} introduced the notion of a fuzzy subset of a
non-empty set $X,$ as a function from $X$ to $[0,1]$. Rosenfeld
\cite{R1} defined the concept of fuzzy group. Since then many papers
have been published in the field of fuzzy algebra.

Recently fuzzy set theory has been well developed in the context of
hyperalgebraic structure theory (for more details see [1-18],
 [20-26],  [29-30], [33-34]).

In \cite{S1} Sen, Ameri and Chowdhury introduced the notions of
fuzzy hypersemigroups and obtained characterization of them. In this
paper we we continuing our previous work \cite{S1}, to  study of
fuzzy $\Gamma$-hypersemigroups. In this context we find some results
of algebraic properties of $\Gamma$-hypersemigroups. In this regards
we introduce $\Gamma$-hyperoperation and study  fuzzy hyperideals.
In particular, we study fuzzy (resp. strongly)regular,  fundamental
relation of $\Gamma$-hypersemigroups.

\section{Preliminaries}

\noindent{\bf Definition 2.1.} Let $M=\{a,b,c,...\}$ and
$\Gamma=\{\alpha,\beta,\gamma,...\}$ be two non-empty sets. Then $M$
is called a $\Gamma$-semigroup if there exists a mapping
$M\times\Gamma\times M\longrightarrow M$ written as
$(a,\gamma,b)\longmapsto a\gamma b$ satisfying the following
identity
\begin{center}
$(a\alpha b)\beta c=a\alpha (b\beta c)$ for all $a,b,c \in M$ and
for all $\alpha,\beta\in \Gamma.$
\end{center}
Let $N$ be a non-empty subset of $M$. Then $N$ is called a sub
$\Gamma$-semigroup of $M$ if $a\gamma b\in N$ for all $a,b \in N$
and $\gamma \in \Gamma$.

\noindent Let $H$ be a non-empty set and let $P^{\ast}(H)$ be the
set of all non-empty subsets of $H$. A hyperoperation on $H$ is a
map $\circ:H\times H \longrightarrow P^{\ast}(H)$ and the couple
$(H,\circ)$ is called a hypergroupoid. If $A$ and $B$ are non-empty
subsets of $H$, then we denote $A\circ B=\bigcup_{a\in A, b\in
B}a\circ b$, $~x\circ A=\{x\}\circ A$ and $ A\circ x=A\circ \{x\}$.

 \noindent {\bf Definition 2.2.} A hypergroupoid $(H,\circ)$ is called a
 semihypergroup if for all $x,y,z $ of $H$ we have $(x\circ y)\circ z=x\circ(y\circ z)$
which means that
\begin{center}
$\bigcup_{u\in x\circ y}~u\circ z=\bigcup_{v\in y\circ z}~x\circ v$
\end{center}

\noindent A semihypergroup $(h,\circ)$ is called a hypergroup if for
all $x\in H$, we have
\begin{center}
$x\circ H=H\circ x=H$.
\end{center}

 \noindent {\bf Definition 2.3.} Let $M$ and $\Gamma$ be two
non-empty sets. $M$ is called a $\Gamma$-semihypergroup if every
$\gamma \in \Gamma$ be a hyperoperation on $M$, i.e, $x\gamma
y\subseteq M$ for every $x,y \in M$, and for every $\alpha,\beta\in
\Gamma$ and $x,y,z \in M$ we have $x\alpha(y\beta z)=(x\alpha
y)\beta z.$

\noindent If every $\gamma \in \Gamma$ is an operation, then $M$ is
a $\Gamma$-semigroup.

\noindent Let $A$ and $B$ be two non-empty subsets of $M$ and
$\gamma \in \Gamma$ we define:
\begin{center}
$A\gamma B=\cup~\{a\gamma b~|~a\in A,b\in B\}.$
\end {center}
 Also
\begin{center}
$A\Gamma B=\cup~\{a\gamma b~|~ a\in A, b\in B ~and~\gamma \in
\Gamma\}=\bigcup_{\gamma \in \Gamma}~A\gamma B.$
\end{center}

\section{Fuzzy $\Gamma$-hyperoperations}

\noindent {\bf Definition 3.1.} Let $M$ and $\Gamma$ be two
non-empty sets. $F(M)$ denote the set of all fuzzy subsets of $M$. A
fuzzy $\Gamma$-hyperoperation on $M$ is a mapping
$\circ:M\times\Gamma\times M\longrightarrow F(M)$ written as
$(a,\gamma,b)\longmapsto a\circ\gamma\circ b$.
 $M$ together with a fuzzy $\Gamma$-hyperoperation  is called a fuzzy
 $\Gamma$-hypergroupoid.

\noindent {\bf Definition 3.2.} A fuzzy $\Gamma$-hypergroupoid
$(M,\circ)$ is called a fuzzy $\Gamma$-hypersemigroup if for all
$a,b,c\in M,~\alpha,\beta\in \Gamma,~(a\circ\alpha\circ
b)\circ\beta\circ c=a\circ\alpha\circ(b\circ\beta\circ c)$ where for
any fuzzy subset $\mu$ of $M$
\begin{center}
$(a\circ\alpha\circ\mu)(r)=\left\{
\begin{array}{cc}
\vee_{t\in M}((a\circ\alpha \circ t)(r)\wedge\mu(t)),&\mu\neq 0\\
0,& otherwise
\end{array}\right.$\\
\end{center}

\noindent and

\begin{center}
      $(\mu\circ\alpha\circ a)(r)=\left\{
     \begin{array}{cc}
      \vee_{t\in M}((\mu(t)\wedge (t\circ\alpha\circ a)(r)),&\mu\neq 0\\
     0,& otherwise
\end{array}\right.$\\
\end{center}

  \noindent{\bf Definition 3.3.}  Let $\mu,\nu$  be two fuzzy subsets of a fuzzy $\Gamma$-hypergroupoid
   $(M,\circ)$, then we define $\mu\circ\gamma\circ\nu$ by $(\mu\circ\gamma\circ\nu)(t)=\vee_{p,q\in M}~(\mu(p)\wedge(p\circ\gamma\circ
   q)(t)\wedge\nu(q))$, for all $t\in M,~\gamma\in \Gamma$.

   \noindent{\bf Definition 3.4.} A fuzzy $\Gamma$-hypersemigroup $(M,\circ)$ is called a fuzzy $\Gamma$-hypergroup if
    $x\circ\gamma \circ M=M\circ \gamma\circ x=\chi_{M}$,
    for all $x\in M,~\gamma \in \Gamma$.

   \noindent {\bf Example 3.5.} Define a fuzzy $\Gamma$-hyperoperation  on a non-empty set
    $M$ by $a\circ \gamma\circ b=\chi_{\{a,\gamma,b\}}$, where $\chi_{\{a,\gamma,b\}}$ denotes the characteristic
    function of the set $\{a,\gamma,b\}$, for
    all $a,b\in M,~\gamma \in \Gamma$, or
    $a\circ \gamma \circ b=\chi_{\{a,b\}}$.
    Then $(M,\circ)$ is a fuzzy $\Gamma$-hypersemigroup.

  \noindent {\bf Example 3.6.}  Let $M$ be a $\Gamma$-semigroup ($\Gamma$-semihypergroup). Define a fuzzy
    $\Gamma$-hyperoperation on $M$ by  $a\circ \gamma\circ b=\chi_{a\gamma b}$, for all $a,b\in M,~\gamma \in \Gamma$,
    where $\chi_{a\gamma b}$ is the
    characteristic function  of the element $M$, then $(M,\circ)$ is a fuzzy
   $\Gamma$-hypersemigroup. If $M$ be a $\Gamma$-group then $(M,\circ)$  is a fuzzy $\Gamma$-hypergroup.

  \noindent {\bf Example 3.7.} Let $M$ be a $\Gamma$-semigroup and $\mu\neq 0$ be a fuzzy $\Gamma$-semigroup on $M$. Let $a,b\in M$ and $\gamma \in\Gamma$.
    Define a fuzzy $\Gamma$-hyperoperation $\circ$ on $M$ by\\
\begin{center}
 $(a\circ \gamma \circ b)(t)=\left\{
\begin{array}{cc}
\mu(a)\wedge\mu(b),~&if~ t=a\gamma b\\
0,& otherwise.
\end{array}\right.$\\
\end{center}

\noindent Then $(M,\circ)$ is a fuzzy $\Gamma$-hypersemigroup.

\noindent {\bf Example 3.8.}  For each positive integer $n$ consider
the set $X_{n}=\{-\infty,0,1,...,n\}$. Define fuzzy
$\Gamma$-hyperoperation on $X_{n}$ and non-empty set $\Gamma$ by
$a\circ\gamma\circ b=\chi_{max\{a,b\}}$ for all $a,b\in X_{n},\gamma
\in \Gamma$ and $-\infty$ is assumed to satisfying  the conditions
that $-\infty\leq a$ for all $a\in X_{n}$ and
$-\infty\circ\gamma\circ a=-\infty$. This gives $X_{n}$ the
structure of a fuzzy $\Gamma$-hypersemigroup.

\noindent {\bf Proof.} Let $a,b,c\in X_{n},\alpha,\beta \in \Gamma$.
$((a\circ\alpha\circ b)\circ\beta\circ c)(t)=(\chi_{max\{a,b\}}\circ
\beta\circ c)(t)=\bigvee_{r\in X_{n}}(\chi_{max\{a,b\}}(r)\wedge
(r\circ\beta\circ c)(t)) $\\

\noindent $=\left\{
\begin{array}{cc}
 (r\circ\beta\circ c)(t),~&r=max\{a,b\}\\
 0, & otherwise.
 \end{array}\right.\\
 =\left\{
 \begin{array}{cc}
 \chi_{max\{r,c\}}(t),&r=max\{a,b\}\\
 0,&otherwise.
 \end{array}\right.\\
 =\left\{
 \begin{array}{cc}
 1,~&t=max\{a,b,c\}\\
 0,& otherwise.
 \end{array}\right.$\\

\noindent  Similarly, we can show that

  \begin{center}
  $(a\circ \alpha\circ( b\circ\beta\circ c))(t)=\left\{
  \begin{array}{cc}
  1,~& t=max\{a,b,c\}\\
  0,& otherwise.
  \end{array}\right.$\\
  \end{center}

  \noindent In this example, if we consider fuzzy
  $\Gamma$-hyperoperation on $X_{n}$ by\\

\begin{center}
  $(a\circ\gamma\circ
  b)(t)=\left\{
  \begin{array}{cc}
  1/2,& t=min\{a+b,n\}\\
  0, & otherwise,
  \end{array}\right.$
  \end{center}
  \noindent where $+$ is integers additive operation, then $(X_{n},\circ)$ is a
fuzzy $\Gamma$-hypersemigroup.$\Box$

\noindent {\bf Example 3.9.} Let $S(M)$ be set of all subsets
nonempty set $M$ and $\Gamma$ be a non-empty set. Define fuzzy
$\Gamma$-hyperoperation $\circ$ on $S(M)$ by

\begin{center}
$(A\circ\gamma\circ B)(C)=\left\{
\begin{array}{cc}
1/3,& C\subseteq A\cup B\\
0, &otherwise.
\end{array}\right.$
\end{center}

\noindent Let $M,\Gamma$ be  nonempty sets, and $M$ endowed with a
fuzzy $\Gamma$-hyperoperation $\circ$ and for all $a,b\in
M,\gamma\in\Gamma$, consider the $p$-cuts $(a\circ\gamma\circ
b)_{p}=\{t\in M:~(a\circ\gamma\circ b)(t)\geq p\}$ of
$a\circ\gamma\circ b$, where $p\in[0,1]$.

\noindent For all $p\in[0,1]$, we define the following crisp
$\Gamma$-hyperoperation on $M$: $a\circ_{p} \gamma \circ_{p}
b=(a\circ \gamma \circ b)_{p}$.

\noindent {\bf Theorem 3.10.} For all $a,b,c,u\in M$ and
$\alpha,\beta\in \Gamma$ and for all $p\in[0,1]$ the following
equivalence holds:

\begin{center}
$(a\circ \alpha \circ (b\circ \beta \circ c))\geq
p\Longleftrightarrow u\in a\circ_{p} \alpha \circ_{p}(b\circ_{p}
\beta \circ_{p} c).$
\end{center}

\noindent {\bf Theorem 3.11.} $(M,\circ)$ is a fuzzy
$\Gamma$-hypersemigroup if and only if $\forall p\in
[0,1],~(M,\circ_{p})$ is a $\Gamma$-hypersemigroup.

\noindent {\bf Proof.} It is obvious.$\Box.$

\noindent {\bf Theorem 3.12.} For all $a\in M$, the following
equivalence holds:
\begin{center}
$a\circ\gamma \circ
M=\chi_{M}\Longleftrightarrow \forall p\in[0,1], ~a\circ _{p}\gamma
\circ_{p} M=M.$
\end{center}

\noindent {\bf Proof.} If  $a\circ\gamma \circ M=\chi_{M}$, then for
all $t\in M$ and  $p\in[0,1]$, we have $\bigvee_{u\in M}(a\circ
\gamma \circ u)(t)=1\geq p$, whence there exists $m\in M$ such that
$(a\circ\gamma\circ m)(t)\geq p$, which means that $t\in
a\circ_{p}\gamma \circ_{p }m$. Hence, $\forall
p\in[0,1],~a\circ_{p}\gamma\circ_{p} M=M$. Conversely, for $p=1$ we
have $a\circ_{1}\gamma\circ_{1} M=M$, whence for all $t\in M$, there
exists $u\in M$, such that $t\in a\circ_{1}\gamma\circ_{1} u$, which
means that $(a\circ\gamma\circ u)(t)=1$. In other words, $a\circ
\gamma \circ M=\chi_{M}$.$\Box$

\noindent {\bf Theorem 3.13.} Let $(M,\circ)$ be a fuzzy
$\Gamma$-hypersemigroup. Then $\chi_{a}\circ \gamma \circ
\chi_{b}=a\circ \gamma \circ b$, for all $a,b\in M, ~\gamma \in
\Gamma$.

\noindent {\bf Proof.} $(\chi_{a}\circ \gamma\circ
\chi_{b})(t)=\vee_{p\in M}~(\chi_{a}(p)\wedge(p\circ \gamma \circ
\chi_{b})(t))=(a\circ\gamma\circ \chi_{b})(t)=\vee_{q\in
M}((a\circ\gamma\circ q)(t)\wedge \chi_{b}(q))=(a\circ\gamma \circ
b)(t)$, for all $t\in M,\gamma \in \Gamma$. Therefore
$\chi_{a}\circ\gamma \circ \chi_{b}=a\circ \gamma \circ b$, for all
$a,b\in M,\gamma \in \Gamma$. Hence the result.$\Box$

\noindent {\bf Theorem 3.14.} Let $(M,\circ)$ be a fuzzy
$\Gamma$-hypersemigroup. Then

\noindent $(i)~~a\circ\alpha\circ
(b\circ\beta\circ\mu)=(a\circ\alpha \circ b)\circ\beta\circ\mu,$ for
all $a,b\in M,\alpha,\beta \in \Gamma$ and for all $\mu \in F(M)$.

\noindent $(ii) ~a\circ\alpha\circ(\mu\circ \beta \circ
b)=(a\circ\alpha\circ\mu)\circ \beta \circ b,$ for all $a,b\in
M,\alpha,\beta \in \Gamma$ and for all $\mu \in F(M)$.

\noindent $(iii)~ \mu\circ \alpha \circ (a\circ\beta\circ
b)=(\mu\circ\alpha\circ a )\circ\beta \circ b$,for all $a,b\in
M,\alpha,\beta \in \Gamma$ and for all $\mu \in F(M)$.

\noindent $(iv)~ \mu \circ\alpha\circ (a
\circ\beta\circ\nu)=(\mu\circ \alpha\circ a)\circ\beta \circ\nu$ for
all $a\in M,\alpha,\beta \in \Gamma$ and for all $\mu,\nu \in F(M)$.

\noindent $(v)~ a\circ\alpha \circ(\mu\circ\beta
\circ\nu)=(a\circ\alpha\circ\mu )\circ\beta\circ\nu$ for all $a\in
M,\alpha,\beta \in \Gamma$ and for all $\mu ,\nu \in F(M)$.

\noindent $(vi)~ \mu\circ\alpha \circ(\nu \circ\beta\circ
a)=(\mu\circ\alpha \circ\nu)\circ\beta \circ a$ for all $a\in
M,\alpha,\beta \in \Gamma$ and for all $\mu,\nu \in F(M)$.

\noindent  $(vii) ~\mu\circ \alpha\circ(\nu\circ\beta
\circ\delta)=(\mu\circ \alpha\circ\nu)\circ\beta\circ\delta$ for all
$\mu,\nu,\delta \in F(M)$.

\noindent {\bf Proof.} It is straight forward.

\section{Fuzzy $\Gamma$-hyperideals}

\noindent {\bf Definition 4.1.} A fuzzy subset $\mu$ of a fuzzy
$\Gamma$-hypersemigroup  $(M,\circ)$ is called a {\it~fuzzy~ sub
$\Gamma$-hypersemigroup} of $(M,\circ)$ if $\mu \circ \gamma
\circ\mu\subseteq \mu$.

\noindent {\bf Theorem 4.2.} If $(M,\circ)$ is a fuzzy
$\Gamma$-hypersemigroup and $\mu,\nu$ are two fuzzy sub
$\Gamma$-hypersemigroups of
 $(M,\circ)$, then $\mu \cap\nu$ is also a fuzzy sub $\Gamma$-hypersemigroup of $(M,\circ)$.

\noindent {\bf Definition 4.3.} A fuzzy subset $\mu$ of a  fuzzy
$\Gamma$-hypersemigroup is called a {\it left fuzzy
$\Gamma$-hyperideal} if $a\circ \gamma \circ \mu\subseteq \mu$, for
all $a\in M,~\gamma \in \Gamma$.

\noindent Similarly we can define a {\it right fuzzy $
\Gamma$-hyperideal} of a fuzzy $\Gamma$-hypersemigroup $(M,\circ)$.

\noindent {\bf Theorem 4.4.} A fuzzy subset $\mu$ of a fuzzy
$\Gamma$-hypersemigroup $(M,\circ)$ is a left fuzzy
$\Gamma$-hyperideal if and only if $M\circ \gamma \circ \mu\subseteq
\mu$.

\noindent {\bf Theorem 4.5.} Let $\mu$ and $\nu$ be two left fuzzy
$\Gamma$-hyperideal of a fuzzy $\Gamma$-hypersemigroup $(M,\circ)$,
then $\mu\cup\nu$ and $\mu\cap\nu$ are also left fuzzy
$\Gamma$-hyperideals of $(M,\circ)$.

\noindent {\bf Theorem 4.6.} Let $(M,\circ)$ be a fuzzy
$\Gamma$-hypersemigroup. Then

\noindent $(i)~~\chi_{M}$ is a left fuzzy $\Gamma$-hyperideal of
$(M,\circ)$, \\
\noindent $(ii)~~\chi_{M}\circ \gamma \circ m=M\circ \gamma \circ
m$, for all
$m\in M,~\gamma \in \Gamma$,\\
\noindent $(iii)~~M\circ \gamma \circ m$ is a left fuzzy
$\Gamma$-hyperideal
of $(M,\circ)$, for all $m\in M,~\gamma \in \Gamma$,\\
\noindent$(iv)$  For any fuzzy subset $\mu\neq 0$ on $M$,$~M\circ
\gamma\circ \mu$ is a left fuzzy $\Gamma$-hyperideal of $(M,\circ)$
for all $\gamma \in \Gamma$.

\noindent {\bf Proof.}
$(i)~m\circ\gamma\circ\chi_{M}\subseteq\chi_{M},$ for all $m\in
M,\gamma \in \Gamma$.

\noindent $(ii)~(\chi_{M}\circ\gamma\circ m)(t)=\vee_{p\in
M}(\chi_{M}(p)\wedge(p\circ\gamma\circ m)(t))=\vee_{p\in
M}(p\circ\gamma\circ m)(t)=(M\circ\gamma \circ m)(t)$, for all $m\in
M,\gamma \in \Gamma$.

\noindent $(iii)~x\circ\alpha\circ(M\circ\gamma\circ
m)=x\circ\alpha\circ (\chi_{M}\circ\gamma\circ m)=(x\circ\alpha\circ
\chi_{M})\circ\gamma\circ m\subseteq \chi_{M}\circ\gamma\circ
m=M\circ\gamma\circ m,$ for all $x,a \in M,\alpha,\gamma \in
\Gamma$.

\noindent
$(iv)~x\circ\alpha\circ(M\circ\gamma\circ\mu)=x\circ\alpha\circ(\chi_{M}\circ\gamma\circ\mu)=
(x\circ\alpha\circ\chi_{M})\circ\gamma\circ\mu\subseteq_{M}\circ\gamma\circ\mu=M
\circ\gamma\circ\mu$, for all $x\in M,\alpha,\gamma \in
\Gamma$.$\Box$

 \noindent {\bf Theorem 4.7.} If $\mu$ is a left
fuzzy $\Gamma$-hyperideal of a fuzzy
 $\Gamma$-suhhypersemigroups $(M,\circ)$, then

 \noindent $(i)~~\mu \circ \gamma \circ m$ is a left fuzzy $\Gamma$-hyperideal of $(M,\circ)$, for all
$m\in M,\gamma \in \Gamma$.

\noindent $(ii)~~\mu \circ \gamma \circ M$ is a left fuzzy
$\Gamma$-hyperideal of $(M,\circ)$.

\noindent {\bf Proof.} $(i)~x\circ\alpha\circ(\mu\circ\gamma\circ
m)=(x\circ\alpha\circ \mu)\circ \gamma\circ m\subseteq
\mu\circ\gamma\circ m$, for all $x\in M,\alpha \in \Gamma$.

\noindent $(ii)~x\circ\alpha\circ(\mu\circ\gamma\circ
M)=x\circ\alpha\circ(\mu\circ\gamma\circ\chi_{M})=(x\circ\alpha\circ\mu)\circ\gamma\circ\chi_{M}\subseteq
\mu\circ\gamma\circ\chi_{M}=\mu\circ\gamma\circ M$, for all $x\in
M,\alpha \in \Gamma$.$\Box$

\noindent {\bf Definition 4.8.} If $\mu\neq0$ is a fuzzy subset of a
fuzzy $\Gamma$-semihypergroup $(M,\circ)$, then the intersection of
all left fuzzy $\Gamma$-hyperideals of $(M,\circ)$ containing $\mu$
($\chi_{M}$ itself being one such is a left fuzzy
$\Gamma$-hyperideal of $(M,\circ)$  containing $\mu$ and contained
in every other such left fuzzy $\Gamma$-hyperideal of $(M,\circ)$.
We call it the left fuzzy $\Gamma$-hyperideal of $(M,\circ)$
generated by $\mu$.

\noindent {\bf Theorem 4.9.} If $\mu\neq0$ is a fuzzy subset of a
fuzzy $\Gamma$-hypersemigroup $(M,\circ)$, then  $\mu\cup(M\circ
\gamma\circ \mu)$ is the smallest left fuzzy $\Gamma$-hyperideal of
$(M,\circ)$ containing $\mu$.

\noindent {\bf Proof.} It is obvious that $\mu\cup(M\circ\gamma\circ
\mu)$ is a left fuzzy $\Gamma$-hyperideal of $(M,\circ)$ containing
$\mu$. Let $\nu$ be a left fuzzy $\Gamma$-hyperideal of $(M,\circ)$
containing $\mu$. Then $\mu\subseteq \nu\Longrightarrow
M\circ\gamma\circ\mu\subseteq
M\circ\gamma\circ\nu\subseteq\nu\Longrightarrow\mu\cup
M\circ\gamma\circ\mu\subseteq\mu\cup\nu=\nu$.

\noindent This shows that $\mu\cup(M\circ \gamma\circ \mu)$ is the
smallest left fuzzy $\Gamma$-hyperideal of $(M,\circ)$ containing
$\mu$.$\Box$.

\noindent{\bf Definition 4.10.} For any fuzzy subset $\mu\neq0$ of a
fuzzy $\Gamma$-hypersemigroup of $(M,\circ)$, $\mu\cup M\circ\gamma
\circ \mu$ is the smallest left fuzzy $\Gamma$-hyperideal of
$(M,\circ)$ containing $\mu$ . It is called the left fuzzy
$\Gamma$-hyperideal of $(M,\circ)$ generated by $\mu$. Similarly we
can define the right fuzzy $\Gamma$-hyperideal of $(M,\circ)$
generated by $\mu$.

\noindent{\bf Definition 4.11.} A fuzzy sub $\Gamma$-hypersemigroup
$\mu$ of a fuzzy $\Gamma$-hypersemigroup $(M,\circ)$ is called a
fuzzy $\Gamma$-hyper bi-ideal of $(M,\circ)$ if $\mu\circ\alpha\circ
y\circ\beta\circ\mu\subseteq\mu $, for all $y\in M,\alpha,\beta\in
\Gamma$.

 \noindent{\bf Theorem 4.12.} If $\mu$ is a fuzzy sub $\Gamma$-hypersemigroup of a fuzzy $\Gamma$-hypersemigroup $(M,\circ)$, then $\mu$ is also a fuzzy
 $\Gamma$-hyper bi-ideal of $(M,\circ)$ if and only if $\mu\circ \alpha \circ M \circ \beta \circ
 \mu\subseteq\mu$, for all $\alpha,\beta\in\Gamma$.

\noindent {\bf Proof.} Let $\mu$ be a fuzzy $\Gamma$-hyper bi-ideal
of $(M,\circ)$. Then $\mu\circ\alpha\circ y
\circ\beta\circ\mu\subseteq \mu$, for all $y\in M$. Now $(\mu
\circ\alpha\circ M\circ \beta\circ\mu)(t)=\vee_{x,y,z\in
M}(x\circ\alpha \circ y\circ \beta\circ
z)(t)\wedge\mu(x)\wedge\mu(z)=\vee_{y\in M}(\mu\circ\alpha\circ
y\circ\beta\circ \mu)(t)\leq \vee_{y \in M}\mu(t)=\mu(t)$, for all
$t\in M,\alpha,\beta \in \Gamma$. Therefore $\mu\circ\alpha\circ
M\circ\beta\circ\mu\subseteq \mu$. Conversely, we suppose
$\mu\circ\alpha\circ M\circ\beta\circ\mu\subseteq \mu$. Then clearly
$\mu \circ\alpha\circ y \circ\beta\circ\mu\subseteq\mu$, for all $y
\in M,\alpha,\beta\in \Gamma$. Hence the result.$\Box.$

 \noindent {\bf Theorem 4.13.} If $\mu$ and $\nu$ are two  fuzzy $\Gamma$-hyper bi-ideals
 of $(M,\circ)$, then $\mu\cap\nu$ is also a  fuzzy $\Gamma$-hyper bi-ideal
 of $(M,\circ)$.

\noindent {\bf Proof.} It is obvious.$\Box$.

 \noindent{\bf Theorem 4.14.} If $\mu$ is a fuzzy subset of a  fuzzy $\Gamma$-hypersemigroup $(M,\circ)$ and  $\nu$ be
 any fuzzy $\Gamma$-hyper bi-ideal
 of $(M,\circ)$, then $\mu\circ\gamma\circ\nu$ and $\nu\circ\gamma\circ\mu$ are both fuzzy $\Gamma$-hyper bi-ideal
 of $(M,\circ)$, for all $\gamma \in \Gamma$.

\noindent {\bf Proof.} $(\mu\circ\gamma\circ\nu)\circ\alpha \circ(
\mu\circ\gamma\circ\nu)=\mu\circ\gamma\circ\nu\circ\alpha\circ\mu\circ\gamma\circ\nu\subseteq
\mu\circ\gamma\circ\nu \circ \alpha \circ
M\circ\gamma\circ\nu\subseteq\mu\circ\gamma\circ\nu$. This implies
that $(\mu\circ\gamma\circ\nu)$ is a fuzzy
$\Gamma$-subhypersemigroup of $(M,\circ)$, for all $\alpha\in
\Gamma$.

\noindent Also, $(\mu\circ\gamma\circ\nu)\circ\alpha\circ
M\circ\beta\circ(\mu\circ\gamma\circ\nu)=\mu\circ\gamma\circ\nu\circ\alpha\circ(M\circ\beta\circ\mu)\circ\gamma\circ\nu
\subseteq\mu\circ\gamma\circ\nu\circ\alpha\circ
M\circ\gamma\circ\nu\subseteq\mu\circ\gamma\circ\nu$. Therefore,
$\mu\circ\gamma\circ\nu$ is a fuzzy $\Gamma$-hyper bi-ideal of
$(M,\circ)$, for all $\alpha,\beta \in\Gamma$.
 \noindent Similarly, we can show that $\nu\circ\gamma\circ\mu$ is a fuzzy $\Gamma$-hyper bi-ideal of
$(M,\circ)$.$\Box$.

 \noindent {\bf Definition 4.15.} A fuzzy subset $\mu$ of a  fuzzy $\Gamma$-hypersemigroup $(M,\circ)$ is called
 $\Gamma$-hyper interior ideal  of $(M,\circ)$ if $x\circ\alpha\circ\mu\circ\beta\circ y \subseteq\mu$, for all $ x,y\in M ,\alpha,\beta \in
 \Gamma$.

\noindent {\bf Theorem 4.16.} A fuzzy subset $\mu$ of a  fuzzy
$\Gamma$-hypersemigroup $(M,\circ)$ is
 $\Gamma$-hyper interior ideal  of $(M,\circ)$ if and only if $M\circ\alpha\circ\mu\circ\beta\circ M\subseteq \mu$, for all $\alpha,\beta \in\Gamma$.

\noindent We can associate a $\Gamma$-hyperoperation on a fuzzy
$\Gamma$-hypersemigroup $(M,\circ)$, as follows:
\begin{center}
$\forall ~a,b \in
M,\gamma\in \Gamma,~a\ast\gamma\ast b=\{x\in M\mid (a\circ\gamma\circ b)(x)>0\}.$\\
\end{center}

\noindent {\bf Theorem 4.17.} If $(M,\circ)$ is a fuzzy
$\Gamma$-hypersemigroup, then $(M,\ast)$ is a
$\Gamma$-hypersemigroup.

\noindent On the other hand, we can define a fuzzy
$\Gamma$-hyperoperation on a $\Gamma$-hypersemigroup $(M,\ast)$, as
follows: $\forall~ a,b\in M,\gamma \in \Gamma,~~a\circ\gamma \circ
b=\chi_{a\ast \gamma\ast b}$.

\noindent {\bf Theorem 4.18.} If $(M,\ast)$ is a
$\Gamma$-hypersemigroup, then $(M,\circ)$ is a fuzzy
$\Gamma$-hypersemigroup.

\noindent Denote by $\mathcal{FHSG}$ the class of all fuzzy
$\Gamma$-hypersemigroups and by $\mathcal{HSG}$ the class of all
$\Gamma$-hypersemigroup. We define
the following two maps:\\
$\varphi:\mathcal{HSG}\longrightarrow\mathcal{FHSG}$,~~$\varphi((M,\ast))=(M,\circ),$
where for all $a,b$ of $M$ and $\gamma \in \Gamma$ we have
$a\circ\gamma\circ b=\chi_{a\ast\gamma\ast b}$ and
$\psi:\mathcal{FHSG}\longrightarrow\mathcal{HSG},~~\psi((M,\circ))=(M,\ast),$
where for all $a,b$ of $M$ and $\gamma \in \Gamma$, we have
$a\ast\gamma\ast b=\{x\mid (a\circ\gamma \circ b)(x)>0\}$.

\noindent {\bf Definition 4.19.} If $\mu_{1},\mu_{2}$ are fuzzy sets
on $M$, then  we say that $\mu_{1}$ is ${\it smaller}$ than
$\mu_{2}$ and we denote $\mu_{1}\leq\mu_{2}$ iff for all $m \in M$,
we have $\mu_{1}(m)\leq\mu_{2}(m)$.

\noindent Let $f:M_{1}\longrightarrow M_{2}$ be a map. If $\mu$ is a
fuzzy set on $M_{1}$, then we define $f(\mu):M_{2}\longrightarrow
[0,1]$, as follows:

\begin{center}
 $(f(\mu))(t)=\bigvee_{r\in f^{-1}(t)}~\mu(r)~~~$
if $~~f^{-1}(t)\neq\phi$,\\
\end{center}

\noindent otherwise we consider $(f(\mu))(t)=0.$

 \noindent{\bf Remark 4.20.}
If $f:M_{1}\longrightarrow M_{2}$ is a map and $m\in M$, then
$f(\chi_{m})=\chi_{f(m)}$. Indeed, for all $t\in M$, we have

\begin{center}
 $(f(\chi_{m}))(t)=\bigvee_{r\in
f^{-1}(t)}~\chi_{m}(r)=\left\{
\begin{array}{cc}
1& ~ f(m)=t\\
0& otherwise.
\end{array}\right.=\chi_{f(m)}(t).$\\
\end{center}

\noindent We can introduce now the fuzzy $\Gamma$-hypersemigroup
homomorphism notion, as follows:

\noindent{\bf Definition 4.21.} Let $(M_{1},\circ_{1})$ and
$(M_{2},\circ_{2})$ be two fuzzy $\Gamma$-hypersemigroups and
$f:M_{1}\longrightarrow M_{2}$ be a map. We say that $f$ is a ${\it
homomorphism}$ of fuzzy $\Gamma$-hypersemigroups if for all $a,b \in
M,~\gamma \in \Gamma$, we have $f(a\circ_{1}\gamma\circ_{1}b)\leq
f(a)\circ_{2}\gamma\circ_{2} f(b).$

\noindent The following two theorems present two connections between
fuzzy $\Gamma$-hypersemigroup homomorphisms and
$\Gamma$-hypersemigroup homomorphism.

\noindent {\bf Theorem 4.22.} Let $(M_{1},\circ_{1})$ and
$(M_{2},\circ_{2})$ be two fuzzy $\Gamma$-hypersemigroups and
$(M_{1},\ast_{1})=\psi(M_{1},\circ_{1}),~(M_{2},\ast_{2})=\psi(M_{2},\circ_{2})$
be the associated $\Gamma$-hypersemigroups. If
$f:M_{1}\longrightarrow M_{2}$ is a homomorphism of fuzzy
$\Gamma$-hypersemigroups, then $f$ is a homomorphism of the
associated $\Gamma$-hypersemigroups, too.

\noindent {\bf Proof}. For all $a,b\in M_{1},\gamma \in \Gamma$, we
have $f(a\circ_{1}\gamma\circ_{1}b)\leq
f(a)\circ_{2}\gamma\circ_{2}f(b)$. Let $x\in
a\ast_{1}\gamma\ast_{1}b$, which means that
$(a\circ_{1}\gamma\circ_{1}b)(x)>0$ and let $t=f(x)$. We have

\begin{center}
$(f(a\circ_{1}\gamma\circ_{1}b))(t)=\bigvee_{r\in f^{-1}(t)}
(a\circ_{1}\gamma\circ_{1}b)(r)\geq
(a\circ_{1}\gamma\circ_{1}b)(x)>0,$
\end{center}
\noindent whence  $(f(a)\circ_{2}\gamma\circ_{2} f(b))(t)>0$. Hence
$t\in f(a)\ast_{2}\gamma\ast_{2}f(b)$. We obtain
$f(a\ast_{1}\gamma\ast_{1}b)\subseteq
f(a)\ast_{2}\gamma\ast_{2}f(b)$.$\Box$

\noindent {\bf Theorem 4.23.} Let $(M_{1},\ast_{1})$ and
$(M_{2},\ast_{2})$ be two $\Gamma$-hypersemigroups and
$(M_{1},\circ_{1})=\varphi(M_{1},\ast_{1}),(M_{2},\circ_{2})=\varphi(M_{2},\ast_{2})$
be the associated fuzzy $\Gamma$-hypersemigroups. The map
$f:M_{1}\longrightarrow M_{2}$ is a homomorphism of
$\Gamma$-hypersemigroups iff it is a homomorphism of fuzzy
$\Gamma$-hypersemigroups.

\noindent {\bf Proof}. $(\Longrightarrow)$ Suppose that $f$ is a
homomorphism of $\Gamma$-hypersemigroups.  Let $a,b\in
M,\gamma\in\Gamma$. For all $t\in Imf$, we have

$(f(a\circ_{1}\gamma\circ_{1} b))(t)=\bigvee_{r\in
f^{-1}(t)}(a\circ_{1}\gamma\circ_{1}b)(r)=\bigvee_{r\in
f^{-1}(t)}\chi_{a\ast_{1}\gamma\ast_{1}b}(r)\\
=\left\{
\begin{array}{cc}
1&~(a\ast_{1}\gamma\ast_{1}b)\cap f^{-1}(t)\neq\phi\\
0&~otherwise\\
\end{array}\right.\\
=\left\{
\begin{array}{cc}
1&~t\in f(a\ast_{1}\gamma\ast_{1}b)\\
0&~otherwise\\
\end{array}\right.$\\
$=\chi_{f(a\ast_{1}\gamma\ast_{1}b)}(t)\leq
\chi_{f(a)\ast_{2}\gamma\ast_{2}f(b)}(t)=(f(a)\circ_{2}\gamma\circ_{2}f(b))(t)$.
If $t\not \in Imf$, then
$(f(a\circ_{1}\gamma\circ_{1}b))(t)=0\leq(f(a)\circ_{2}\gamma\circ_{2}f(b))(t)$.
Hence, $f(a\circ_{2}\gamma\circ_{2}b) \leq
f(a)\circ_{2}\gamma\circ_{2}f(b)$.

\noindent Hence, $f$ is a homomorphism of fuzzy
$\Gamma$-hypersemigroups.\\
 $(\Longleftarrow)$ Conversely, suppose
that $f$ is a homomorphism of fuzzy $\Gamma$-hypersemigroups and
$a,b\in M_{1}$. Then, for all $t\in M_{1}$,we have

\begin{center}
$(f(a\circ_{1}\gamma\circ_{1}b))(t)\leq(f(a)\circ_{2}\gamma\circ_{2}f(b))(t)$,
 \end{center}

\noindent whence we obtain

\begin{center}
$\chi_{f(a\ast_{1}\gamma\ast_{1}b)}(t)\leq
\chi_{f(a)\ast_{2}\gamma\ast_{2}f(b)}(t),$
 \end{center}

\noindent which means that

 \begin{center}
$f(a\ast_{1}\gamma\ast_{1}b)\subseteq
f(a)\ast_{2}\gamma\ast_{2}f(b).$
\end{center}

\noindent Hence $f$ is a homomorphism of $\Gamma$-hypersemigroups.

\section{Fuzzy  fuzzy (strongly) regular relations}

\noindent In \cite{A12} fuzzy regular relations are introduced in
the context of fuzzy hypersemigroups. We define these relations on a
fuzzy $\Gamma$-hypersemigroup:

\noindent Let $\rho$ be an equivalence relation on a fuzzy
$\Gamma$-hypersemigroup $(M,\circ)$ and let $\mu,\nu$ be two fuzzy
subsets on $M$. We say that $\mu\rho\nu$ if the following two
conditions hold:\\
$(1)$ if $\mu(a)>0$, then there exists $b\in M$, such that
$\nu(b)>0$ and $a\rho b$;

\noindent $(2)$ if $\nu(x)>0$, then there exists $y\in M$, such that
$\mu(y)>0$ and $x\rho y$.

\noindent An equivalence relation $\rho$ on a fuzzy
$\Gamma$-hypersemigroup $(M,\circ)$ is called a {\it fuzzy
 $\Gamma$-regular~ relation}  (or a {\it fuzzy
 $\Gamma$-hypercongruence})  on $(M,\circ)$ if, for all $a,b,c \in M,\gamma \in \Gamma$, the
 following implication holds:

 \begin{center}
 $a\rho b\Longrightarrow (a\circ\gamma\circ c)~\rho~(b\circ\gamma\circ
 c)$ and $(c\circ\gamma\circ a)~\rho~ (c\circ\gamma \circ b)$.
\end{center}

 \noindent This condition
is equivalent to

\noindent $a\rho a^{\prime}, b\rho b^{\prime}$ implies
$(a\circ\gamma\circ b)\rho(a^{\prime}\circ\gamma\circ b^{\prime})$
for all $a,b,a^{\prime},b^{\prime}$ of $M$ and $\gamma \in \Gamma$.

\noindent Let $(M,\circ)$ be a fuzzy $\Gamma$-hypersemigroup and let
$\psi(M,\circ)=(M,\ast)$ be the associated $\Gamma$-hypersemigroup,
where, for all $a,b\in M,\gamma \in \Gamma$, we have $a\ast
\gamma\ast b=\{ x\in M~|~(a\circ\gamma\circ b)(x)>0\}$.

\noindent {\bf Theorem 5.1.} An equivalence relation $\rho$ is a
fuzzy $\Gamma$-regular relation on $(M,\circ)$ if and only if $\rho$
is a $\Gamma$-regular relation on $(M,\ast)$.

\noindent {\bf Definition 5.2.} An equivalence relation $\rho$ on a
fuzzy $\Gamma$-hypersemigroup $(M,\circ)$ is called a {\it fuzzy
$\Gamma$-strongly regular relation} on $(M,\circ)$ if, for all
$a,a^{\prime},b,b^{\prime}$ of $M$ and for all $\gamma \in \Gamma$,
such that $a\rho b$ and $a^{\prime}\rho b^{\prime}$, the following
condition holds: \\
$\forall x\in M$ such that $(a\circ\gamma\circ c)(x)>0$ and $\forall
y\in M$ such that $(b\circ\gamma\circ d)(y)>0$, we have $x\rho y$.

\noindent Notice that if $\rho$ is a fuzzy $\Gamma$-strongly
relation on a fuzzy $\Gamma$-hypersemigroup $(M,\circ)$, then it is
fuzzy $\Gamma$-regular on $(M,\circ)$.

\noindent{\bf Theorem 5.3.} An equivalence relation $\rho$ is a
fuzzy $\Gamma$-strongly regular relation on $(M,\circ)$ if and only
if $\rho$ is a $\Gamma$-strongly regular relation on $(M,\ast)$.

\noindent Let $(M,\circ)$ be a fuzzy $\Gamma$-hypersemigroup. Let
$\rho$ be a fuzzy $\Gamma$-hypercongruence on $M$. We define a
$\Gamma$-hyperoperation $\ast$ on $M$ by

\begin{center}
$a\ast \gamma \ast b=\{x\in M:(a\circ \gamma \circ b)(x)>0\}$.
\end{center}

\noindent Therefore, $(M,\ast)$ is a $\Gamma$-hypersemigroup and
$\rho$ is a $\Gamma$-hypercongruence on $(M,\ast)$.

\noindent Let $M/\rho=\{a\rho :a\in M\}$. We define a
$\Gamma$-hyperoperation $\otimes$ on $M/\rho$ by

\begin{center}
 $a\rho\otimes \gamma \otimes b\rho=\{c\rho~:c\in a\ast \gamma
\ast b\}=\{c \rho: (a\circ\gamma\circ b)(c)>0\}$,
\end{center}
\noindent then $(M/\rho,\otimes)$ is a $\Gamma$-hypersemigroup.

\noindent The $\Gamma$-hypersemigroup $(M/\rho,\otimes)$ is called
the quotient $\Gamma$-hypersemigroup induced by the
$\Gamma$-hypercongruence $\rho$ on $(M,\ast)$.

\noindent {\bf Theorem 5.4.} If $\rho$ is a fuzzy $\Gamma$-strong
hypercongruence on a fuzzy $\Gamma$-hypersemigroup $(M,\circ)$, then
$M/\rho=\{a\rho :a\in M\}$ is a fuzzy $\Gamma$-hypersemigroup.

\noindent {\bf Proof.} We define a fuzzy $\Gamma$-hyperoperation
$\ast$ on $M/\rho$ by $(a\rho\ast \gamma \ast
b\rho)(c\rho)=\vee_{a^{\prime}\in a \rho,b^{\prime}\in
b\rho,c^{\prime}\in c\rho}~(a^{\prime}\circ\gamma\circ
b^{\prime})(c^{\prime})$, for all $a\rho, b\rho, c\rho \in
M/\rho,\gamma \in \Gamma$. Clearly $\ast$ is well-defined.

\noindent Now $((a\rho\ast \gamma \ast b\rho)\ast\gamma\ast c\rho)(
d\rho)=\vee_{p\rho\in M/\rho}~((a\rho\ast \gamma \ast
b\rho)(p\rho)\wedge(p\rho\ast \gamma \ast c\rho)(d\rho))=\vee_{p\rho
\in M/\rho}(\vee_{a^{\prime}\in a\rho,b^{\prime}\in
b\rho,p^{\prime}\in p\rho}(a^{\prime}~\circ ~\gamma~ \circ~
b^{\prime})(p^{\prime})\wedge(\vee_{p_{1}\in p\rho, c^{\prime}\in
c\rho,d^{\prime}\in d\rho}~(p_{1}~\circ~ \gamma ~\circ~
c^{\prime})(d^{\prime}))=\vee_{p\in M}(~\vee_{a^{\prime}\in a\rho,~
b^{\prime}\in b\rho,
 ~c^{\prime}\in c\rho,
 ~d_{1}\in d\rho,
p^{\prime}\in p\rho}~((a^{\prime}~\circ~\gamma~\circ
~b^{\prime})(p^{\prime})~\wedge~(p^{\prime}~\circ~\gamma~\circ
~c^{\prime})(d_{1})))=
  \vee_{a^{\prime}\in a\rho,~b^{\prime}\in
b\rho,~c^{\prime}\in c \rho,~d_{1}\in
d\rho}~((a^{\prime}\circ\gamma\circ b^{\prime})\circ\gamma\circ
c^{\prime})(d^{\prime}))=\vee_{a^{\prime}\in a\rho,b^{\prime}\in
b\rho,c^{\prime}\in c \rho,d_{1}\in
d\rho}(a^{\prime}\circ\gamma\circ b^{\prime}\circ \gamma \circ
 c^{\prime})(d_{1}),$

\noindent for all $d\rho\in M/\rho$.

 \noindent Similarly,
we show that $(a\rho\ast\gamma\ast(b\rho\ast\gamma\ast
c\rho))(d\rho)=\vee_{a^{\prime}\in a\rho,b^{\prime}\in
b\rho,c^{\prime}\in c \rho,d_{1}\in
d\rho}(a^{\prime}\circ\gamma\circ b^{\prime}\circ\gamma\circ
c^{\prime})(d_{1})$, for all $d\rho\in M/\rho$.

\noindent Therefore $(a\rho\ast\gamma\ast
b\rho)\ast\gamma\ast(c\rho)=a\rho\ast\gamma\ast(b\rho\ast\gamma\ast
c\rho)$. Therefore, $(M/\rho,\ast)$
 is a fuzzy
$\Gamma$-hypersemigroup.$\Box$

\noindent {\bf Theorem 5.5.} Let $(M,\circ)$ be a fuzzy
$\Gamma$-hypersemigroup and $\psi(M,\circ)=(M,\ast)$ be the
associated $\Gamma$-hypersemigroup. Then we have:

\noindent $(i)$ The relation $\rho$ is a fuzzy $\Gamma$-regular
relation on $(M,\circ)$ if and only if $(M/\rho,\otimes)$ is a
$\Gamma$-hypersemigroup.

\noindent $(ii)$ The relation $\rho$ is a fuzzy $\Gamma$-strongly
regular relation on $(M,\circ)$ if and only if $(M/\rho,\otimes)$ is
a $\Gamma$-semigroup.

\noindent {\bf Proof.} Straightforward.$\Box$

\noindent\textbf{Acknowledgements.}

\noindent\emph{ The first author partially has been supported by the
"Research Center in Algebraic Hyperstructures and Fuzzy Mathematics,
University of Mazandaran, Babolsar, Iran" and "Algebraic
Hyperstructure Excellence, Tarbiat Modares University, Tehran,
Iran".} \centerline{}

\end{document}